\documentclass[12pt]{article}
\usepackage{graphicx}
\usepackage{hyperref,url}
\usepackage{amsthm}
\usepackage{amsmath,amssymb}
\usepackage{mathrsfs}
\usepackage{amsthm}
\usepackage{bibentry}
\usepackage{enumitem}
\usepackage[nottoc, notlof, notlot]{tocbibind}
\usepackage{xcolor}
\usepackage{textcomp}
\numberwithin{equation}{section}
\newtheorem{theorem}{Theorem}[section]
\newtheorem{prop}[theorem]{Proposition}
\newtheorem{lem}[theorem]{Lemma}

\newtheorem{defi}[theorem]{Definition}

\theoremstyle{remark}
\newtheorem{rem}[theorem]{Remark}
\newcommand{\email}[1]{\href{mailto:#1}{\textsf{#1}}}

\hypersetup{pdfborder={0 0 0}}

\newcommand{\R}{\mathbb{R}}

\newcommand{\N}{\mathbb{N}}
\newcommand{\C}{\mathbb{C}}
\begin{document}
\date{}
\theoremstyle{plain}

\title{Global weak solution of 3D-NSE with exponential damping}
\maketitleÂ²
%----------classification, keywords, date

\begin{center}
Jamel BENAMEUR  \\ King Saud university\\
This Project was funded by the National Plan for Science, Technology and Innovation (MAARIFAH), King Abdulaziz City for Science and Technology, Kingdom of Saudi Arabia, Award Number (14-MAT730-02)\\
\email{\sl jamelbenameur@gmail.com}
\end{center}

\begin{abstract}
In this paper we prove the global existence of incompressible Navier-Stokes equations with damping $\alpha (e^{\beta |u|^2}-1)u$, where we use Friedrich method and some new tools. The delicate problem in the construction of a global solution, is the passage to the limit in exponential nonlinear term. To solve this problem, we use a polynomial approximation of the damping part and a new type of interpolation between $L^\infty(\R^+,L^2(\R^3))$ and the space of functions $f$ such that $(e^{\beta|f|^2}-1)|f|^2\in L^1(\R^3)$. Fourier analysis and standard techniques are used.
\end{abstract}

MSC class: 35-XX, 35Q30, 76N10\\
keywords: Navier-Stokes Equations; Friedrich method; weak solution

\maketitle

\section{\bf Introduction}
In this paper we study the global existence of weak solution to the modified incompressible Navier-Stokes equations in three spatial dimensions
$$(S)
  \begin{cases}
     \partial_t u
 -\nu\Delta u+ u.\nabla u  +\alpha (e^{\beta |u|^2}-1)u =\;\;-\nabla p\hbox{ in } \mathbb R^+\times \mathbb R^3\\
     {\rm div}\, u = 0 \hbox{ in } \mathbb R^+\times \mathbb R^3\\
    u(0,x) =u^0(x) \;\;\hbox{ in }\mathbb R^3,
  \end{cases}
$$
where $u=u(t,x)=(u_1,u_2,u_3)$, $p=p(t,x)$ denote respectively the unknown velocity and the unknown pressure of the fluid at the point $(t,x)\in \mathbb R^+\times \mathbb R^3$, and the viscosity of fluid $\nu>0$. Also $\alpha,\beta>0$ denote the parameters of damping term. The terms $(u.\nabla u):=u_1\partial_1 u+u_2\partial_2 u+u_3\partial_3u$, while $u^0=(u_1^o(x),u_2^o(x),u_3^o(x))$ is an initial given velocity. If $u^0$ is quite regular, the divergence free condition determines the pressure $p$. We recall in our case it was assumed the viscosity is unitary ($\nu=1$) in order to simplify the calculations and the proofs of our results. The global existence of weak solution of initial value problem of classical incompressible Navier-Stokes were proved by Leray and Hopf (see \cite{6}-\cite{11}) long before. Uniqueness remains an open problem for the dimensions $d\geq3$. The damping is from the resistance to the motion of the flow. It describes various physical situations such as porous media flow, drag or friction effects, and some dissipative mechanisms (see \cite{1,2,7,8} and references therein). The polynomial damping $\alpha|u|^{\beta-1}u$ is studied in \cite{01} by Cai and Jiu, where they proved the globale existence of weak solution in
 $$L^\infty(\R^+,L^2(\R^3))\cap L^2(\R^+,\dot H^1(\R^3))\cap L^{\beta+1}(\R^+,L^{\beta+1}(\R^3)).$$

 The purpose of this paper is to study the well-posedness of the incompressible Navier-Stokes equations with damping $\alpha (e^{\beta|u|^2}-1)u$.
We will show that the Cauchy problem $(S)$ has global weak solutions for any $\alpha,\beta\in(0,\infty)$. We apply the Friedrich method to construct the approximate solutions and make more delicate estimates
to proceed to compactness arguments. In particular, we obtain new more a priori estimates,
$$\|u(t)\|_{L^2}^2+2\int_0^t\|\nabla u(z)\|_{L^2}^2dz+2\alpha\int_0^t\|(e^{\beta |u(z)|^2}-1)|u(z)|^2\|_{L^1}dz\leq \|u^0\|_{L^2}^2,$$
comparing with the
Navier-Stokes equations, to guarantee that the solution $u$ belongs to
$$L^\infty(\R^+,L^2(\R^3))\cap L^2(\R^+,\dot H^1(\R^3))\cap \mathcal E_\beta$$
where
$$\mathcal E_\beta=\{f:\R^+\times\R^3\rightarrow\R\;{\rm measurable};\;(e^{\beta|f|^2}-1)|f|^2\in L^1(\R^+\times\R^3)\}.$$
Before treating the global existence, we give the definition of solution of $(S)$.
\begin{defi} The function pair $(u,p)$ is called a weak solution of the problem $(S)$ if for any $T>0$, the
following conditions are satisfied:
\begin{enumerate}
\item $u\in L^\infty([0,T] ;L^2_\sigma(\R^3))\cap L^2([0,T];\dot H^1(\R^3))$ and $(e^{\beta|u|^2}-1)|u|^2\in L^1([0,T],L^1(\R^3))$.
\item $\partial_t u
 -\Delta u+ u.\nabla u  +\alpha (e^{\beta |u|^2}-1)u =-\nabla p$ in $\mathcal D'([0,T]\times\R^3)$: for any $\Phi\in C_0^\infty([0,T]\times\R^3)$ such that ${\rm div}\,\Phi(t,x)=0$ for all $(t,x)\in[0,T]\times\R^3$ and $\Phi(T)=0$, we have
$$-\int_0^T(u;\partial_t\Phi)+\int_0^T(\nabla u;\nabla \Phi)+\int_0^T((u.\nabla)u;\Phi)+\alpha\int_0^T((e^{\beta|u|^2}-1)u;\Phi)=(u^0;\Phi(0)).$$
\item ${\rm div}\, u(x, t) = 0$ for a.e. $(t,x)\in[0,T]\times\R^3.$\\\\
Here $(\; ;\; )$ means the inner product in $L^2(\R^3)$.
\end{enumerate}
\end{defi}
\begin{rem} {\bf a)} If $(u,p)$ is a solution of $(S)$, then $(e^{\beta|u|^2}-1)u\in L^1_{loc}(\R^+,L^1(\R^3))$. Indeed: For $T>0$, put the following subsets of $[0,T]\times\R^3$
$$\begin{array}{lcl}
A_{T}^1&=&\{(t,x)\in[0,T]\times\R^3;\;0<|u(t,x)|< 1\},\\
A_{T}^2&=&\{(t,x)\in[0,T]\times\R^3;\;|u(t,x)|\geq1\}.
\end{array}$$
Then
$$\begin{array}{l}
\displaystyle\int_{[0,T]\times\R^3}(e^{\beta|u|^2}-1)|u|=\displaystyle\int_{A_{T}^1\cup A_{T}^2}(e^{\beta|u|^2}-1)|u|\\
\leq\displaystyle\int_{A_{T}^1}(e^{\beta|u|^2}-1)|u|+\int_{A_{T}^2}(e^{\beta|u|^2}-1)|u|\\
\leq\displaystyle\int_{A_{T}^1}\frac{(e^{\beta|u|^2}-1)}{|u|}|u|^2+\int_{A_{T}^2}(e^{\beta|u|^2}-1)\frac{1}{|u|}|u|^2\\
\leq\displaystyle M_\beta\int_{A_{T}^1}|u|^2+\int_{A_{T}^2}(e^{\beta|u|^2}-1)|u|^2,\;M_\beta=\max_{0<X\leq 1}\frac{(e^{\beta X^2}-1)}{X}\\
\leq\displaystyle M_\beta\int_0^T\int_{\R^3}|u|^2+\int_0^T\int_{\R^3}(e^{\beta|u|^2}-1)|u|^2\\
\leq\displaystyle M_\beta T\|u\|_{L^\infty([0,T],L^2(\R^3))}+\int_{[0,T]\times\R^3}(e^{\beta|u|^2}-1)|u|^2\\
\leq\displaystyle M_\beta T\|u\|_{L^\infty([0,T],L^2(\R^3))}+\|(e^{\beta|u|^2}-1)|u|^2\|_{L^1([0,T],L^1(\R^3))},
\end{array}$$
which implies the desired result.\\
{\bf b)} By this definition the scalar pressure function $p$ is well defined in $L^1_{loc}(\R^+,H^{-s}(\R^3))$, for all $s>3/2$. Indeed: Formally, if $p$ exists, then
$$p=(-\Delta)^{-1}\Big({\rm div}\,u.\nabla u  +\alpha {\rm div}\,(e^{\beta |u|^2}-1)u\Big).$$
So, just prove that $(-\Delta)^{-1}({\rm div}\,(u.\nabla u)),  (-\Delta)^{-1}(\alpha {\rm div}\,(e^{\beta |u|^2}-1)u)\in L^1_{loc}(\R^+,H^{-s}(\R^3))$, for all $s>3/2$. For this, for $s>3/2$ we have
$$\begin{array}{l}
\|(-\Delta)^{-1}({\rm div}\,(u.\nabla u))(t)\|_{H^{-s}(\R^3)}=\|(-\Delta)^{-1}({\rm div}\,({\rm div}(u\otimes u)))(t)\|_{H^{-s}(\R^3)}\\
\leq\|u\otimes u(t)\|_{H^{-s}(\R^3)}\\
\leq\displaystyle\Big(\int_{\R^3} (1+|\xi|^2)^{-s}|\mathcal F(u\otimes u)(t,\xi)|^2d\xi\Big)^{1/2}\\
\leq\displaystyle c_s\|\mathcal F(u\otimes u)(t)\|_{L^\infty(\R^3)},\;c_s=\Big(\int_{\R^3} (1+|\xi|^2)^{-s}d\xi\Big)^{1/2}\\
\leq c_s\|u\otimes u(t)\|_{L^1(\R^3)}\\
\leq c_s\|u(t)\|_{L^2(\R^3)}^2\in L^1_{loc}(\R^+).
\end{array}$$
In other hand, we have
$$\begin{array}{l}
\|(-\Delta)^{-1}({\rm div}\,(e^{\beta|u|^2}-1)u(t)\|_{H^{-s}(\R^3)}\\
\leq\displaystyle\Big(\int_{\R^3}(1+|\xi|^2)^{-s}|\xi|^{-2}|\mathcal F((e^{\beta|u|^2}-1)u)(t,\xi)|^2d\xi\Big)^{1/2}\\
\leq\displaystyle C_s\|\mathcal F((e^{\beta|u|^2}-1)u(t))\|_{L^\infty(\R^3)},\;C_s=\Big(\int_{\R^3} (1+|\xi|^2)^{-s}|\xi|^{-2}d\xi\Big)^{1/2}\\
\leq C_s\|(e^{\beta|u|^2}-1)u(t)\|_{L^1(\R^3)}\in L^1_{loc}(\R^+).
\end{array}$$
Then $p$ is well defined in $L^1_{loc}(\R^+,H^{-s}(\R^3)),\,s>3/2$.\\
{\bf c)} If $(u(t,x),p(t,x))$ is a solution of $(S)$, then $(u(t,x),p(t,x)+h(t))$ is also solution  of $(S)$. Then, if we look for $p$ in the space $L^1_{loc}(\R^+,H^{-2}(\R^3))$, we get $h=0$ and we can express $p$ in terms of $u$.
\end{rem}
In our case of exponential damping, we are trying to find more regularity of Leray solution in $\cap_pL^p(\R^+,L^p(\R^3))$. In particular, we give a new energy estimate. Our main result is the following.
\begin{theorem}\label{th1} Let $u^0\in L^2(\mathbb R^3)$ be a divergence free vector fields, then there is a global solution of $(S)$
$u\in L^\infty(\R^+,L^2(\mathbb R^3)\cap C(\R^+,H^{-2}(\R^3))\cap L^2(\R^+,\dot H^1(\mathbb R^3))\cap\mathcal E_\beta$. Moreover, for all $t\geq0$ \begin{equation}\label{eqth1}\|u(t)\|_{L^2}^2+2\int_0^t\|\nabla u(z)\|_{L^2}^2dz+2\alpha\int_0^t\|(e^{\beta |u(z)|^2}-1)|u(z)|^2\|_{L^1}dz\leq \|u^0\|_{L^2}^2.\end{equation}
\end{theorem}
\begin{rem} {\bf a)} If $(u,p)$ is a solution of $(S)$ given by Theorem \ref{th1}, then by (\ref{eqth1}) and Proposition \ref{prop1}-(3), we get
$$\lim_{t\rightarrow0}\|u(t)-u^0\|_{L^2}=0.$$
{\bf b)} The fact $(e^{\beta |u|^2}-1)|u|^2\in L^1(\R^+,L^1(\R^3))$ implies that $u\in \cap _{4\leq p<\infty}L^p(\R^+,L^p(\R^3))$. Indeed: we have $$\int_0^\infty\|(e^{\beta|u(t)|^2}-1)|u(t)|^2\|_{L^1}dt=\sum_{k=4}^{\infty}\frac{\beta^k}{k!}\int_0^\infty\|u(t)\|_{L^{2k+2}}^{2k+2}dt.$$
{\bf c)} The fact $u\in L^\infty(\R^+,L^2(\mathbb R^3)\cap C(\R^+,H^{-2}(\R^3))$ implies that $u\in C_b(\R^+,H^{-r}(\R^3))$, for all $r>0$.
\end{rem}
The remainder of our paper is organized as follows. In the second section we give some notations, definitions and preliminary results. Section 3 is devoted to prove Theorem \ref{th1}, this proof is done in two steps. In the first, we give a general result of equicontinuity, in the second we apply Friedritch method to construct a global solution of $(S)$, the only delicate point is to show the convergence of the nonlinear part $(e^{\beta|u_{\varphi(n)}|^2}-1)u_{\varphi(n)}$ to $(e^{\beta|u|^2}-1)u$. To do this, we approximate $(e^{\beta|u_{\varphi(n)}|^2}-1)u_{\varphi(n)}$ by
$\displaystyle\sum_{k=1}^m\frac{\beta^k}{k!}|u_{\varphi(n)}|^{2k}u_{\varphi(n)}$ and we use the convergence in $L^p_{loc}(\R^+,L^p_{loc}(\R^3))$ given by technical lemmas.
\section{\bf Notations and preliminary results}
\subsection{Notations} In this section, we collect some notations and definitions that will be used later.\\
\begin{enumerate}
\item[$\bullet$] For $R>0$: $B_R=\{x\in\R^3;\;|x|<R\}$.
\item[$\bullet$] For $m\in\N$: $\displaystyle P_m(X)=\sum_{k=1}^m\frac{X^k}{k!}$. Clearly, for each $R>0$, we have $P_m(\beta R^2)\nearrow (e^{\beta R^2}-1)$, as $m\rightarrow+\infty$.
\item[$\bullet$] The Fourier transformation is normalized as
$$
\mathcal{F}(f)(\xi)=\widehat{f}(\xi)=\int_{\mathbb R^3}\exp(-ix.\xi)f(x)dx,\,\,\,\xi=(\xi_1,\xi_2,\xi_3)\in\mathbb R^3.
$$
\item[$\bullet$] The inverse Fourier formula is
$$
\mathcal{F}^{-1}(g)(x)=(2\pi)^{-3}\int_{\mathbb R^3}\exp(i\xi.x)g(\xi)d\xi,\,\,\,x=(x_1,x_2,x_3)\in\mathbb R^3.
$$
\item[$\bullet$] The convolution product of a suitable pair of function $f$ and $g$ on $\mathbb R^3$ is given by
$$
(f\ast g)(x):=\int_{\mathbb R^3}f(y)g(x-y)dy.
$$
\item[$\bullet$] If $f=(f_1,f_2,f_3)$ and $g=(g_1,g_2,g_3)$ are two vector fields, we set
$$
f\otimes g:=(g_1f,g_2f,g_3f),
$$
and
$$
{\rm div}\,(f\otimes g):=({\rm div}\,(g_1f),{\rm div}\,(g_2f),{\rm div}\,(g_3f)).
$$
Moreover, if $\rm{div}\,g=0$ we obtain
$$
{\rm div}\,(f\otimes g):=g_1\partial_1f+g_2\partial_2f+g_3\partial_3f:=g.\nabla f.
$$
\item[$\bullet$] Let $(B,||.||)$, be a Banach space, $1\leq p \leq\infty$ and  $T>0$. We define $L^p_T(B)$ the space of all
measurable functions $[0,t]\ni t\mapsto f(t) \in B$ such that $t\mapsto||f(t)||\in L^p([0,T])$.\\
\item[$\bullet$] The Sobolev space $H^s(\mathbb R^3)=\{f\in \mathcal S'(\mathbb R^3);\;(1+|\xi|^2)^{s/2}\widehat{f}\in L^2(\mathbb R^3)\}$.\\
\item[$\bullet$] The homogeneous Sobolev space $\dot H^s(\mathbb R^3)=\{f\in \mathcal S'(\mathbb R^3);\;\widehat{f}\in L^1_{loc}\;{\rm and}\;|\xi|^s\widehat{f}\in L^2(\mathbb R^3)\}$.
\item[$\bullet$] $L^2_\sigma(\R^3)=\{f\in (L^2(\R^3))^3;\;{\rm div}\,f=0\}$.
\item[$\bullet$] For $R>0$, the Friedritch operator $J_R$ is defined by
$$J_R(D)f=\mathcal F^{-1}({\bf 1}_{B_R}(\xi)\widehat{f}).$$
\item[$\bullet$] The Leray projector $\mathbb P:(L^2(\R^3))^3\rightarrow (L^2(\R^3))^3$ is defined by
$$\mathcal F(\mathbb P f)=\widehat{f}(\xi)-(\widehat{f}(\xi).\frac{\xi}{|\xi|})\frac{\xi}{|\xi|}=M(\xi)\widehat{f}(\xi);\;M(\xi)=(\delta_{k,l}-\frac{\xi_k\xi_l}{|\xi|^2})_{1\leq k,l\leq 3}.$$
\item[$\bullet$] For $R>0$, we define the operator $A_R(D)$ by
$$A_R(D)f=\mathbb P J_R(D)f=\mathcal F^{-1}(M(\xi){\bf 1}_{B_R}(\xi)\widehat{f}).$$
\item[$\bullet$] If $p,q\in (1,\infty)$ and $\Omega_1\subset\R^n,\;\Omega_2\subset\R^m$ are two open subsets, we define $L^p(\Omega_1,L^q(\Omega_2))$ by:
$f\in L^p(\Omega_1,L^q(\Omega_2))$ if $f:\Omega_1\times\Omega_2\rightarrow\C,\;(x_1,x_2)\rightarrow f(x_1,x_2)$\; is measurable and
$$\|f\|_{L^p(\Omega_1,L^q(\Omega_2))}=\|\|f\|_{L^q(\Omega_2)}\|_{L^p(\Omega_1)}$$ is well defined and finite.
Particularly, if $p=q$, then $L^p(\Omega_1,L^q(\Omega_2))=L^p(\Omega_1\times\Omega_2)$ and $$\|f\|_{L^p(\Omega_1,L^p(\Omega_2))}=\|f\|_{L^p(\Omega_1\times\Omega_2)}.$$
\end{enumerate}
\subsection{Preliminary results}
In this section, we recall some classical results and we give new technical lemmas.
\begin{prop}(\cite{HBAF})\label{prop1} Let $H$ be Hilbert space.
\begin{enumerate}
\item If $(x_n)$ is a bounded sequence of elements in $H$, then there is a subsequence $(x_{\varphi(n)})$ such that
$$(x_{\varphi(n)}|y)\rightarrow (x|y),\;\forall y\in H.$$
\item If $x\in H$ and $(x_n)$ is a bounded sequence of elements in $H$ such that
$$(x_n|y)\rightarrow (x|y),\;\forall y\in H,$$
then $\|x\|\leq\displaystyle\liminf_{n\rightarrow\infty}\|x_n\|.$
\item If $x\in H$ and $(x_n)$ is a bounded sequence of elements in $H$ such that
$$(x_n|y)\rightarrow (x|y),\;\forall y\in H$$
and
$$\limsup_{n\rightarrow\infty}\|x_n\|\leq \|x\|,$$
then $\displaystyle\lim_{n\rightarrow\infty}\|x_n-x\|=0.$
\end{enumerate}
\end{prop}
\begin{lem}(\cite{JYC})\label{LP}
Let $s_1,\ s_2$ be two real numbers and $d\in\N$.
\begin{enumerate}
\item If $s_1<d/2$\; and\; $s_1+s_2>0$, there exists a constant  $C_1=C_1(d,s_1,s_2)$, such that: if $f,g\in \dot{H}^{s_1}(\mathbb{R}^d)\cap \dot{H}^{s_2}(\mathbb{R}^d)$, then $f.g \in \dot{H}^{s_1+s_2-\frac{d}{2}}(\mathbb{R}^d)$ and
$$\|fg\|_{\dot{H}^{s_1+s_2-\frac{d}{2}}}\leq C_1 (\|f\|_{\dot{H}^{s_1}}\|g\|_{\dot{H}^{s_2}}+\|f\|_{\dot{H}^{s_2}}\|g\|_{\dot{H}^{s_1}}).$$
\item If $s_1,s_2<d/2$\; and\; $s_1+s_2>0$ there exists a constant $C_2=C_2(d,s_1,s_2)$ such that: if $f \in \dot{H}^{s_1}(\mathbb{R}^d)$\; and\; $g\in\dot{H}^{s_2}(\mathbb{R}^d)$, then  $f.g \in \dot{H}^{s_1+s_2-\frac{d}{2}}(\mathbb{R}^d)$ and
$$\|fg\|_{\dot{H}^{s_1+s_2-\frac{d}{2}}}\leq C_2 \|f\|_{\dot{H}^{s_1}}\|g\|_{\dot{H}^{s_2}}.$$
\end{enumerate}
 \end{lem}
\begin{lem}\label{lempn1} For $m\in\N$ and $\beta>0$, there is $c_{m,\beta}>0$ such that
$$\Big(P_m(\beta z^2)\Big)^2\leq c_{m,\beta}(e^{\beta z^2}-1)z^2,\;\forall z\in\R.$$
\end{lem}
{\bf Proof.} It suffices to prove that $\Big(P_m(\beta z^2)\Big)^2(e^{\beta z^2}-1)^{-1}z^{-2}$ is bounded on $(0,\infty)$.
\begin{lem}\label{lempn2} For $m\in\N$ and $\beta>0$, there is $C_{m,\beta}>0$ such that
$$|P_m(\beta |x|^2)x-P_m(\beta |y|^2)y|\leq C_{m,\beta}\Big(P_m(\beta |x|^2))+P_m(\beta |y|^2)\Big)|x-y|,\;\forall x,y\in\R^3.$$
\end{lem}
{\bf Proof.} Suppose that $|y|\leq |x|$. We have
$$P_m(\beta |x|^2)x-P_m(\beta |y|^2)y=P_m(\beta |x|^2)(x-y)+\sum_{k=1}^m\frac{\beta^k}{k!}(|x|^{2k}-|y|^{2k})y.$$
For $k\in\{1,...,m\}$, we have
$$\begin{array}{lcl}\Big||x|^{2k}-|y|^{2k}\Big|.|y|&\leq& 2k|x-y|.(|x|^{2k-1}+|y|^{2k-1}).|y|\\
&\leq& 2m|x-y|.(|x|^{2k-1}|y|+|y|^{2k})\\
&\leq& 2m|x-y|.(|x|^{2k}+|y|^{2k}),
\end{array}$$
which implies the desired result.
\begin{lem}\label{lff1} For all $s>d/2$: $L^1(\R^d)\hookrightarrow H^{-s}(\R^d)$. Moreover, we have
$$\|f\|_{H^{-s}(\R^d)}\leq \sigma_{s,d}\|f\|_{L^1(\R^d)},\;\forall f\in H^{-s}(\R^d),$$
where $\sigma_{s,d}=\Big(\int_{\R^d}(1+|\xi|^2)^{-s}d\xi\Big)^{1/2}$.
\end{lem}
 {\bf Proof.} Using the fact $s>d/2$, we get
 $$\int_{\R^d}(1+|\xi|^2)^{-s}|\widehat{f}(\xi)|^2d\xi\leq\displaystyle\Big(\int_{\R^d}(1+|\xi|^2)^{-s}d\xi\Big)\|\widehat{f}\|_{L^\infty(\R^d)}^2
 \leq\displaystyle\Big(\int_{\R^d}(1+|\xi|^2)^{-s}d\xi\Big)\|f\|_{L^1(\R^d)}^2,$$
which ends the proof.
\begin{lem}\label{l122} Let $\varepsilon\in(0,1)$ and $d\in\N$ such that $d\geq2$, then there is a constant $C>0$ such that: If $(f,g)\in H^{-\frac{\varepsilon}{2}}(\R^d)\times H^{1-\frac{\varepsilon}{2}}(\R^d)$, then $fg\in H^{1-\varepsilon-\frac{d}{2}}(\R^d)$ and
$$\|fg\|_{H^{1-\varepsilon-\frac{d}{2}}(\R^d)}\leq C\|f\|_{H^{-\frac{\varepsilon}{2}}(\R^d)}\|g\|_{H^{1-\frac{\varepsilon}{2}}(\R^d)}.$$
\end{lem}
{\bf Proof.} We have
$$\begin{array}{lll}
\|fg\|_{H^{1-\varepsilon-\frac{d}{2}}}^2&\leq&\|fg\|_{\dot H^{1-\varepsilon-\frac{d}{2}}}^2\\
&\leq&\displaystyle\int_{\xi}|\xi|^{2(1-\varepsilon-\frac{d}{2})}\Big|\int_{\eta}\widehat{f}(\eta)\widehat{g}(\xi-\eta)\Big|^2\\
&\leq&\displaystyle\int_{\xi}|\xi|^{2(1-\varepsilon-\frac{d}{2})}
\Big(\int_{\eta}|\widehat{f}(\xi-\eta)||\widehat{g}(\eta)|\Big)^2\\
&\leq&\displaystyle\int_{\xi}|\xi|^{2(1-\varepsilon-\frac{d}{2})}
\Big(\int_{\eta}(1+|\xi-\eta|^2)^{-\varepsilon/4}|\widehat{f}(\xi-\eta)|(1+|\xi-\eta|^2)^{\varepsilon/4}|\widehat{g}(\eta)|\Big)^2.
\end{array}$$
By using the elementary inequality
$$(1+|\xi-\eta|^2)^{\frac{\varepsilon}{4}}\leq (1+2|\xi|^2+2|\eta|^2)^{\frac{\varepsilon}{4}}\leq 6^{\frac{\varepsilon}{4}}(1+|\xi|^{\frac{\varepsilon}{2}}+|\eta|^{\frac{\varepsilon}{2}}),$$
we get
$$\begin{array}{lll}
\|fg\|_{H^{1-\varepsilon-\frac{d}{2}}}^2&\leq&\displaystyle C\int_{\xi}|\xi|^{2(1-\varepsilon-\frac{d}{2})}
\Big(\int_{\eta}(1+|\xi-\eta|^2)^{-\frac{\varepsilon}{4}}|\widehat{f}(\xi-\eta)|(1+|\xi|^{\frac{\varepsilon}{2}}+|\eta|^{\varepsilon/2})|\widehat{g}(\eta)|\Big)^2\\
&\leq&C(I_1+I_2+I_3),
\end{array}$$
with
$$\begin{array}{lll}
I_1&=&\displaystyle\int_{\xi}|\xi|^{2(1-\varepsilon-\frac{d}{2})}
\Big(\int_{\eta}(1+|\xi-\eta|^2)^{-\frac{\varepsilon}{4}}|\widehat{f}(\xi-\eta)||\widehat{g}(\eta)|\Big)^2=\|u_1v_1\|_{\dot H^{1-\varepsilon-\frac{d}{2}}}^2\\
u_1&=&\mathcal F^{-1}((1+|\xi|^2)^{-\frac{\varepsilon}{4}}|\widehat{f}(\xi)|)\\
v_1&=&\mathcal F^{-1}(|\widehat{g}(\xi)|),
\end{array}$$
$$\begin{array}{lll}
I_2&=&\displaystyle\int_{\xi}|\xi|^{2(1-\frac{\varepsilon}{2}-\frac{d}{2})}
\Big(\int_{\eta}(1+|\xi-\eta|^2)^{-\frac{\varepsilon}{4}}|\widehat{f}(\xi-\eta)||\widehat{g}(\eta)|d\eta\Big)^2d\xi=\|u_2v_2\|_{\dot H^{1-\frac{\varepsilon}{2}-\frac{d}{2}}}^2\\
u_2&=&\mathcal F^{-1}((1+|\xi|^2)^{-\frac{\varepsilon}{4}}|\widehat{f}(\xi)|)\\
v_2&=&\mathcal F^{-1}(|\widehat{g}(\xi)|)
\end{array}$$
and
$$\begin{array}{lll}
I_3&=&\displaystyle\int_{\xi}|\xi|^{2(1-\varepsilon-\frac{d}{2})}
\Big(\int_{\eta}(1+|\xi-\eta|^2)^{-\frac{\varepsilon}{4}}|\widehat{f}(\xi-\eta)|
|\eta|^{\frac{\varepsilon}{2}}|\widehat{g}(\eta)|\Big)^2=\|u_3v_3\|_{\dot H^{1-\varepsilon-\frac{d}{2}}}^2\\
u_3&=&\mathcal F^{-1}((1+|\xi|^2)^{-\frac{\varepsilon}{4}}\widehat{f}(\xi))\\
v_3&=&\mathcal F^{-1}(|\xi|^{\frac{\varepsilon}{2}}|\widehat{g}(\xi)|).
\end{array}$$
By applying Lemma \ref{LP} with the following choices
$$(s_1,t_1)=(0,1-\varepsilon),\;(s_2,t_2)=(0,1-\frac{\varepsilon}{2}),\;(s_3,t_3)=(0,1-\varepsilon)$$
we get
$$\|u_iv_i\|_{\dot H^{s_i+t_i-\frac{d}{2}}(\R^d)}\leq C\|u_i\|_{\dot H^{s_i}(\R^d)}\|v_i\|_{\dot H^{t_i}(\R^d)},\;i\in\{1,2,3\}.$$
Then
$$\begin{array}{l}
I_1\leq C\|u_1\|_{\dot H^0}^2\|v_1\|_{\dot H^{1-\varepsilon}}^2\leq C\|f\|_{H^{-\frac{\varepsilon}{2}}}^2\|g\|_{\dot H^{1-\varepsilon}}^2\leq C\|f\|_{H^{-\frac{\varepsilon}{2}}}^2\|g\|_{H^{1-\varepsilon}}^2\leq C\|f\|_{H^{-\frac{\varepsilon}{2}}}^2\|g\|_{H^{1-\frac{\varepsilon}{2}}}^2\\\\
I_2\leq C\|u_2\|_{\dot H^0}^2\|v_2\|_{\dot H^{1-\frac{\varepsilon}{2}}}^2\leq C\|f\|_{H^{-\frac{\varepsilon}{2}}}^2\|g\|_{H^{1-\frac{\varepsilon}{2}}}^2\\\\
I_3\leq C\|u_3\|_{\dot H^0}^2\|v_3\|_{\dot H^{1-\varepsilon}}^2\leq C\|f\|_{H^{-\frac{\varepsilon}{2}}}^2\|g\|_{\dot H^{1-\frac{\varepsilon}{2}}}^2\leq C\|f\|_{H^{-\frac{\varepsilon}{2}}}^2\|g\|_{H^{1-\frac{\varepsilon}{2}}}^2,
\end{array}$$
which ends the proof of the Lemma.
\begin{lem}\label{lem11} Let $p\in(1,+\infty)$ and $\Omega\neq\emptyset$ be an open subset of $\R^4$. If $(f_n)$ is a bounded sequence in $L^2(\Omega)\cap L^p(\Omega)$, then there is a subsequence $(f_{h(n)})$ and $f\in L^2(\Omega)\cap L^p(\Omega)$ such that
$$\begin{array}{rcl}
\displaystyle\lim_{n\rightarrow\infty}(f_{h(n)}|g)_{L^2}&=&(f|g)_{L^2},\;\;\;\forall g\in L^2(\Omega)\\
\|f\|_{L^p}&\leq&\displaystyle\liminf_{n\rightarrow\infty}\|f_{h(n)}\|_{L^p}\\
\|f\|_{L^2}&\leq&\displaystyle\liminf_{n\rightarrow\infty}\|f_{h(n)}\|_{L^2}.
\end{array}$$
\end{lem}
\noindent{\bf Proof.} $(f_n)$ is a bounded sequence in the Hilbert space $L^2(\Omega)$, then there is a subsequence $(f_{h(n)})$ and $f\in L^2(\Omega)$ such that
$$\lim_{n\rightarrow\infty}(f_{h(n)}|g)_{L^2}=(f|g)_{L^2},\;g\in L^2(\Omega).$$
Particularly, for $g\in C_0^\infty(\Omega)$, we have
$$|(f_{h(n)}|g)_{L^2}|\leq \|f_{h(n)}\|_{L^p}\|g\|_{L^q},\;\frac{1}{p}+\frac{1}{q}=1.$$
Then
$$\liminf_{n\rightarrow\infty}|(f_{h(n)}|g)_{L^2}|\leq (\liminf_{n\rightarrow\infty}\|f_{h(n)}\|_{L^p})\|g\|_{L^q}$$and$$|(f|g)_{L^2}|\leq (\liminf_{n\rightarrow\infty}\|f_{h(n)}\|_{L^p})\|g\|_{L^q}.$$
As $ C_0^\infty(\Omega)$ is dense in $L^q(\Omega)$, then $f\in L^p(\Omega)$ and $\|f\|_{L^p}\leq \displaystyle\liminf_{n\rightarrow\infty}\|f_{n}\|_{L^p}$.\\
The last result is giving by Proposition \ref{prop1}.
\begin{rem}\label{rem11} Let $p\in(1,+\infty)$ and $\Omega\neq\emptyset$ be an open subset of $\R^4$. If $(f_n)$ is a bounded sequence in $L^2(\Omega)\cap L^p(\Omega)$ and $f\in L^2(\Omega)$ such that
$$\lim_{n\rightarrow\infty}(f_n|g)_{L^2}=(f|g)_{L^2},\;\;\;\forall g\in L^2(\Omega).$$
Then $f\in L^p(\Omega)$ and
$$\begin{array}{lcl}
\|f\|_{L^2}&\leq&\displaystyle\liminf_{n\rightarrow\infty}\|f_{n}\|_{L^2}\\
\|f\|_{L^p}&\leq&\displaystyle\liminf_{n\rightarrow\infty}\|f_{n}\|_{L^p}.
\end{array}$$
\end{rem}
\begin{lem}\label{lem12} Let $T>0$ and $(f_n)$ be a bounded sequence in $L^2_T(H^1(\R^3))$ such that
 $$f_n\rightarrow f\;strongly\;in\;C([0,T],H^{-4}_{loc}(\R^3)).$$
Then $$f_n\longrightarrow f\;strongly\,in\;L^2([0,T],L^2_{loc}(\R^3)).$$
\end{lem}
\noindent{\bf Proof.} Combining the inclusion $C([0,T],H^{-4}_{loc}(\R^3))\hookrightarrow L^2([0,T],H^{-4}_{loc}(\R^3))$ and the interpolation $$L^2([0,T],H^{-4}_{loc}(\R^3)\cap H^1(\R^3))\hookrightarrow L^2([0,T],L^2_{loc}(\R^3)),$$ we obtain the desired result.
\begin{lem}\label{lem2.9} Let $T>0$ and $p_0\in(2,+\infty)$. If $(f_n)$ is a bounded sequence in $L^2_T(H^1(\R^3))\cap L^{p_0}_T(L^{p_0}(\R^3))$ such that
$$f_n\rightarrow f\;strongly\;in\;C([0,T],H^{-4}_{loc}(\R^3)),$$
then $f\in L^{p_0}_T(L^{p_0}(\R^3))$ and
$$f_n\longrightarrow f\; strongly\,in\;L^{p}([0,T],L^p_{loc}(\R^3)),\;\forall 2\leq p<p_0.$$
\end{lem}
\noindent{\bf Proof.} By Lemma \ref{lem12}, we get
$$f_n\longrightarrow f\;{\rm strongly\,in}\;L^{2}([0,T],L^2_{loc}(\R^3)).$$
By interpolation between $L^{p_0}([0,T],L^{p_0}(\R^3))$ and $L^{2}([0,T],L^2_{loc}(\R^3))$, we get
$$f_n\rightarrow f\;{\rm strongly\;in}\;L^{p}([0,T],L^{p}_{loc}(\R^3)),\;\forall p\in[2,p_0).$$
It remains to show that $f\in L^{p_0}_T(L^{p_0}(\R^3))$. For $m\in\N$, put
$$\Omega_m=]0,T[\times B_m\;{\rm and}\;g_{m,n}={\bf 1}_{\Omega_m}f_n$$
Clearly $\Omega_m$ is an open subset of $\R^4$ and the sequence $(g_{m,n})_{n}$ is bounded in
$$L^2([0,T],L^2(B_m))\cap L^{p_0}([0,T],L^{p_0}(B_m))=L^2(\Omega_m)\cap L^{p_0}(\Omega_m).$$
Then, by Lemma \ref{lem11} and Remark \ref{rem11}, we get ${\bf 1}_{\Omega_m}f\in L^2(\Omega_m)\cap L^{p_0}(\Omega_m)$ and
$$\|{\bf 1}_{\Omega_m}f\|_{L^{p_0}(\Omega_m)}\leq \liminf_{n\rightarrow\infty}\|g_{n,m}\|_{L^{p_0}(\Omega_m)}.$$
But $\|g_{n,m}\|_{L^{p_0}(\Omega_m)}=\|f_n\|_{L^{p_0}(\Omega_m)}\leq \|f_n\|_{L^{p_0}_T(L^{p_0}(\R^3))}$, then
$$\|{\bf 1}_{\Omega_m}f\|_{L^{p_0}(\Omega_m)}\leq \liminf_{n\rightarrow\infty}\|f_n\|_{L^{p_0}_T(L^{p_0}(\R^3))}.$$
By applying the Monotonic Convergence Theorem, we get $f\in L^{p_0}_T(L^{p_0}(\R^3))$ and
$$\|f\|_{L^{p_0}_T(L^{p_0}(\R^3))}\leq \liminf_{n\rightarrow\infty}\|f_n\|_{L^{p_0}_T(L^{p_0}(\R^3))}.$$
\section{\bf Proof of Theorem \ref{th1}.}
This proof is done in two steps:
\subsection{Step 1} In this step we prove a general result for bounded sequence in energy space of the system $(S)$.
\begin{prop}\label{prop2} Let $\nu_1,\nu_2,\nu_3\in[0,\infty)$, $r_1,r_2,r_3\in(0,\infty)$ and $f_0\in L^2_\sigma(\R^3)$. For $n\in\N$, let $F_n:\R^+\times\R^3\rightarrow\R^3$ be a measurable function in $C^1(\R^+,L^2(\R^3))$ such that $$A_n(D)F_n=F_n,\;F_n(0,x) =A_n(D)f_0(x)$$ and
$$\begin{array}{l}
\displaystyle(E1)\;\;\;\;\partial_t F_n+\sum_{k=1}^3\nu_k|D_k|^{2r_k} F_n+ A_n(D){\rm div}\,(F_n\otimes F_n)+\alpha A_n(D)[(e^{\beta |F_n|^2}-1)F_n] =0\\
\displaystyle(E2)\;\;\;\;\displaystyle\|F_n(t)\|_{L^2}^2+2\sum_{k=1}^3\nu_k\int_0^t\||D_k|^{r_k} F_n\|_{L^2}^2+2\alpha\int_0^t\|(e^{\beta |F_n|^2}-1)|F_n|^2\|_{L^1}\leq \|f_0\|_{L^2}^2.
 \end{array}$$
Then:
\begin{enumerate}
\item[$\bullet$] For $m\in\N$, there is a constant $C=C(m,\alpha,\beta)$ such that
\begin{equation}\label{ps1eq1}\int_0^\infty\int_{\R^3}\Big(P_m(\beta|F_n(t,x)|^2)\Big)^2dtdx\leq C\|f_0\|_{L^2}^2.\end{equation}
\item[$\bullet$] For all $T,R>0$ and $s>3/2$, we have
\begin{equation}\label{ps1eq11}\int_0^T\|(e^{\beta |F_n|^2}-1)|F_n|\|_{H^{-s}}\leq \sigma_{3,s}(M_{\beta,R}T+\frac{1}{2R\alpha})\|f_0\|_{L^2}^2,\;\;\forall n\in\N,\end{equation}
where $\displaystyle M_{\beta,R}=\sup_{0<r\leq R}\frac{e^{\beta r^2}-1}{r}$.
\item[$\bullet$] For every $\varepsilon>0$ there is $\delta=\delta(\varepsilon,\alpha,\beta,\nu_1,\nu_2,\nu_3,r_1,r_2,r_3,\|f_0\|_{L^2})>0$ such that: For all $t_1,t_2\in\R^+$, we have
\begin{equation}\label{ps1eq2}\Big(|t_2-t_1|<\delta\Longrightarrow \|F_n(t_2)-F_n(t_1)\|_{H^{-s_0}}<\varepsilon\Big),\;\forall n\in\N,\end{equation}
with $\displaystyle s_0=\max(3,2\max_{1\leq i\leq 3}r_i).$
\end{enumerate}
\end{prop}
\begin{rem} As $\|F_n(t)\|_{L^2}\leq \|f_0\|_{L^2}$ for all $(t,n)\in\R^+\times\N$, then by interpolation between $H^{-s_0}(\R^3)$ and $H^0(\R^3)=L^2(\R^3)$ we get: for all $\sigma>0$, for each $\varepsilon>0$, there is $\delta>0$ such that
$$\Big(|t_2-t_1|<\delta\Longrightarrow \|F_n(t_2)-F_n(t_1)\|_{H^{-\sigma}}<\varepsilon\Big),\;\forall n\in\N.$$
\end{rem}
\noindent{\bf Proof of Proposition \ref{prop2}.} $\bullet$ Firstly, (\ref{ps1eq1}) is given by $(E2)$ and Lemma \ref{lempn1}.\\
$\bullet$ To prove (\ref{ps1eq11}), beginning by noting that the function $h_\beta:[0,+\infty)\longrightarrow\R^+$, $y\longmapsto\left\{\begin{array}{l}
\displaystyle\frac{e^{\beta y^2}-1}{y}\;{\rm if}\,t>0\\
0\;{\rm if}\,t=0\end{array}\right.$ is continuous and $M_{\beta,R}=\sup_{y\in[0,R]}h_\beta(y)$ is well defined and finite.\\
For $t\in[0,T]$, put the following subset of $\R^3$
$$X_{n,R,t}=\{x\in\R^3;\;|F_n(t,x)|\leq R\}.$$
Then, for $x\in\R^3$ we have
$$\begin{array}{lcl}
x\in X_{n,R,t}&\Longrightarrow& (e^{\beta |F_n(t,x)|^2}-1)|F_n(t,x)|=h_\beta(|F_n(z,x)|)|F_n(t,x)|^2\leq M_{\beta,R}|F_n(t,x)|^2\\\\
x\in X_{n,R,t}^c&\Longrightarrow& \displaystyle(e^{\beta |F_n(t,x)|^2}-1)|F_n(t,x)|\leq \frac{1}{R}(e^{\beta |F_n(t,x)|^2}-1)|F_n(t,x)|^2.
\end{array}$$
By Lemma \ref{lff1} we obtain
$$
D_n(T):=\displaystyle\int_{0}^{T}\|(e^{\beta |F_n(z)|^2}-1)|F_n(z)|\|_{H^{-s}}\leq\displaystyle \sigma_{3,s}\int_{0}^{T}\|(e^{\beta |F_n(z)|^2}-1)|F_n(z)|\|_{L^1}dz.$$
Combining this inequality with the above inequalities, we get
$$\begin{array}{lcl}
\displaystyle
\frac{1}{\sigma_{3,s}}D_n(T)&=&\displaystyle\int_{0}^{T}\int_{\R^3}(e^{\beta |F_n(z,x)|^2}-1)|F_n(z,x)|\\
&=&\displaystyle\int_{0}^{T}\int_{X_{n,R,z}}(e^{\beta |F_n(z,x)|^2}-1)|F_n(z,x)|+\int_{0}^{T}\int_{X_{n,R,z}^c}(e^{\beta |F_n(z,x)|^2}-1)|F_n(z,x)|\\
&\leq&\displaystyle M_{\beta,R}\int_{0}^{T}\int_{X_{n,R,z}}|F_n(z,x)|^2+\frac{1}{R}\int_{0}^{T}\int_{X_{n,R,z}^c}(e^{\beta |F_n(z,x)|^2}-1)|F_n(z,x)|^2\\
&\leq&\displaystyle M_{\beta,R}\int_{0}^{T}\int_{\R^3}|F_n(z,x)|^2+\frac{1}{R}\int_{0}^{T}\int_{\R^3}(e^{\beta |F_n(z,x)|^2}-1)|F_n(z,x)|^2\\
&\leq&\displaystyle M_{\beta,R}\int_{0}^{T}\|F_n(z)\|_{L^2}^2dz+\frac{1}{R}\int_{0}^{T}\|(e^{\beta |F_n(z)|^2}-1)|F_n(z)|^2\|_{L^1(\R^3)}.
\end{array}$$
By using $(E_2)$, we get the desired result.\\
$\bullet$ To prove (\ref{ps1eq2}), integrate $(E1)$ over $[t_1,t_2]\subset\R^+$, we get
$$\|F_n(t_2)-F_n(t_1)\|_{H^{-s_0}}\leq I_{n,1}(t_1,t_2)+I_{n,2}(t_1,t_2)+I_{n,3}(t_1,t_2),$$
with
$$\begin{array}{lcl}
I_{n,1}(t_1,t_2)&=&\displaystyle\sum_{k=1}^3\nu_k\int_{t_1}^{t_2}\||D_k|^{2r_k}F_n\|_{H^{-s_0}}\\
I_{n,2}(t_1,t_2)&=&\displaystyle\int_{t_1}^{t_2}\|A_n(D){\rm div}\,(F_n\otimes F_n)\|_{H^{-s_0}}\\
I_{n,3}(t_1,t_2)&=&\displaystyle\int_{t_1}^{t_2}\|\alpha A_n(D)[(e^{\beta |F_n|^2}-1)F_n]\|_{H^{-s_0}}.
\end{array}$$
Let $\varepsilon>0$ be a positive real, let us find a positive real $\delta>0$ such that if $|t_2-t_1|<\delta$, we get
$$I_{n,k}(t_1,t_2)<\varepsilon/3,\;\;k\in\{1,2,3\}.$$
$\bullet$ To estimate $I_{n,1}(t_1,t_2)$, we write
$$I_{n,1}(t_1,t_2)\leq \sum_{k=1}^3\nu_k\int_{t_1}^{t_2}\|F_n\|_{H^{2r_k-s_0}}
\leq (\sum_{k=1}^3\nu_k)\int_{t_1}^{t_2}\|F_n(z)\|_{H^0}dz\leq (\sum_{k=1}^3\nu_k)\|\widehat{f}_0\|_{L^2}(t_2-t_1).$$
Then if $|t_2-t_1|<\delta_1:=\displaystyle\frac{\varepsilon}{3(\displaystyle\sum_{k=1}^3\nu_k)\|\widehat{f}_0\|_{L^2}+3}$, we get $I_{n,1}(t_1,t_2)<\varepsilon/3$.\\
$\bullet$ Estimate $I_{n,2}(t_1,t_2)$: By Lemma \ref{lff1}, we get
$$\begin{array}{lcl}
\displaystyle I_{n,2}(t_1,t_2)&\leq&\displaystyle\int_{t_1}^{t_2}\|{\rm div}\,(F_n\otimes F_n)(z)\|_{H^{-3}}dz\\
&\leq&\displaystyle\int_{t_1}^{t_2}\|(F_n\otimes F_n)(z)\|_{H^{-2}}dz\\
&\leq&\displaystyle \sigma_{2,3}\int_{t_1}^{t_2}\|(F_n\otimes F_n)(z)\|_{L^1}dz\\
&\leq&\displaystyle \sigma_{2,3}\int_{t_1}^{t_2}\|F_n(z)\|_{L^2}^2dz\\
&\leq&\displaystyle \sigma_{2,3}\|f_0\|_{L^2}^2(t_2-t_1).
\end{array}$$
Then if $|t_2-t_1|<\delta_2:=\displaystyle\frac{\varepsilon}{3\sigma_{2,3}\|f_0\|_{L^2}^2+3}$, we get $I_{n,2}(t_1,t_2)<\varepsilon/3$.\\
\noindent$\bullet$ Estimate of $I_{n,3}(t_1,t_2)$: Let $R>0$ (to fixed later) and $t_1<t_2\in\R^+$. By Lemma \ref{lff1} we obtain
$$
\displaystyle\int_{t_1}^{t_2}\|(e^{\beta |F_n(z)|^2}-1)|F_n(z)|\|_{H^{-3}}\leq\displaystyle \sigma_{3,3}\int_{t_1}^{t_2}\|(e^{\beta |f(z)|^2}-1)|f(z)|\|_{L^1}dz.$$
Combining this inequality with the above inequalities, we get
$$\begin{array}{lcl}
\displaystyle
\frac{I_3(t_1,t_2)}{\alpha \sigma_{3,3}}&=&\displaystyle\int_{t_1}^{t_2}\int_{\R^3}(e^{\beta |F_n(z,x)|^2}-1)|F_n(z,x)|\\
&=&\displaystyle\int_{t_1}^{t_2}\int_{X_{n,R,z}}(e^{\beta |F_n(z,x)|^2}-1)|F_n(z,x)|+\int_{t_1}^{t_2}\int_{X_{n,R,z}^c}(e^{\beta |F_n(z,x)|^2}-1)|F_n(z,x)|\\
&\leq&\displaystyle M_{\beta,R}\int_{t_1}^{t_2}\int_{X_{n,R,z}}|F_n(z,x)|^2+\frac{1}{R}\int_{t_1}^{t_2}\int_{X_{n,R,z}^c}(e^{\beta |F_n(z,x)|^2}-1)|F_n(z,x)|^2\\
&\leq&\displaystyle M_{\beta,R}\int_{t_1}^{t_2}\int_{\R^3}|F_n(z,x)|^2+\frac{1}{R}\int_{t_1}^{t_2}\int_{\R^3}(e^{\beta |F_n(z,x)|^2}-1)|F_n(z,x)|^2\\
&\leq&\displaystyle M_{\beta,R}\int_{t_1}^{t_2}\|F_n(z)\|_{L^2}^2dz+\frac{1}{R}\int_{t_1}^{t_2}\|(e^{\beta |F_n(z)|^2}-1)|F_n(z)|^2\|_{L^1(\R^3)}.
\end{array}$$
By using $(E_2)$, we get
\begin{equation}\label{eqtr}\frac{I_3(t_1,t_2)}{\alpha \sigma_{3,3}}\leq\displaystyle M_{\beta,R}\|f_0\|_{L^2}^2(t_2-t_1)+\frac{1}{2R\alpha}\|f_0\|_{L^2}^2.\end{equation}
Then, with the choices $R=R_\varepsilon=\displaystyle\frac{3\sigma_{3,3}\|f_0\|_{L^2}^2+3}{\varepsilon}$ and $\displaystyle \delta_3=\frac{\varepsilon}{6\alpha \sigma_{3,3} M_{\beta,R_\varepsilon}\|f_0\|_{L^2}^2+6}$, we get
$$|t_1-t_2|<\delta_3\Longrightarrow I_{n,3}(t_1,t_2)<\varepsilon/3.$$
To conclude, it suffices to take $\delta=\min\{\delta_1,\delta_2,\delta_3\}.$

\subsection{Step 2} In this step we construct a global solution of $(S)$, where we use a method inspired by \cite{JYC}. For this, consider the approximate system with the parameter $n\in\N$:
$$(S_n)
  \begin{cases}
     \partial_t u
 -\Delta J_nu+ J_n(J_nu.\nabla J_nu  +\alpha J_n[(e^{\beta |J_nu|^2}-1)J_nu] =\;\;-\nabla p_n\hbox{ in } \mathbb R^+\times \mathbb R^3\\
 p_n=(-\Delta)^{-1}\Big({\rm div}\,J_n(J_nu.\nabla J_nu  +\alpha {\rm div}\,J_n[(e^{\beta |J_nu|^2}-1)J_nu]\Big)\\
     {\rm div}\, u = 0 \hbox{ in } \mathbb R^+\times \mathbb R^3\\
    u(0,x) =J_nu^0(x) \;\;\hbox{ in }\mathbb R^3.
  \end{cases}
$$
\begin{enumerate}
\item[$\bullet$] By Cauchy-Lipschitz theorem, we obtain a unique solution $u_n\in C^1(\R^+,L^2_\sigma(\R^3))$ of $(S_n)$. Moreover, $J_nu_n=u_n$ and
\begin{equation}\label{theq1}\|u_n(t)\|_{L^2}^2+2\int_0^t\|\nabla u_n\|_{L^2}^2+2\alpha\int_0^t\|(e^{\beta |u_n|^2}-1)|u_n|^2\|_{L^1}\leq \|u^0\|_{L^2}^2.\end{equation}
\item[$\bullet$] By inequality (\ref{theq1}), we get $(u_n)$ is bounded in $L^2_{loc}(\R^+,H^1(\R^3))$:
\begin{equation}\label{theq2}\int_0^T(\|u_n\|_{L^2}^2+\|\nabla u_n\|_{L^2}^2)\leq \|u^0\|_{L^2}^2T+\frac{\|u^0\|_{L^2}^2}{2}=M_T,\;\;\forall T\geq0.\end{equation}
\item[$\bullet$] Using inequality (\ref{theq1}) and the fact
\begin{equation}\label{theq3}(e^{\beta |u_n|^2}-1)|u_n|^2=\sum_{k=1}^\infty\frac{\beta^k}{k!}|u_n|^{2k+2}\end{equation}
we get: for all $k\in\N$, $(u_n)$ is bounded in $L^{2k+2}(\R^+\times\R^3)$ and
\begin{equation}\label{theq4}\int_0^\infty\int_{\R^3}|u_n|^{2k+2}\leq \frac{k!2\alpha}{\beta^k}.\end{equation}
\item[$\bullet$] By applying Proposition \ref{prop2} and Remark \ref{rem11} we get
\begin{equation}\label{theq5} {\rm the\,sequence}\;(u_n)\;{\rm is\, equicontinuous\,in}\;C_b(\R^+,H^{-1}(\R^3)).\end{equation}
\item[$\bullet$] Let $(T_q)_q\in(0,\infty)^\N$ such that $T_q<T_{q+1}$ and $T_q\rightarrow\infty$ as $q\rightarrow\infty$. Let $(\theta_q)_{q\in\N}$ be a sequence in $C_0^\infty(\R^3)$ such that: for all $q\in\N$
$$\left\{\begin{array}{l}
\theta_q(x)=1,\;\forall x\in B(0,q+1+\frac{1}{4})\\
\theta_q(x)=0,\;\forall x\in B(0,q+2)^c\\
0\leq \theta_q\leq 1.
\end{array}\right.$$
Using (\ref{theq1})-(\ref{theq5}) and classical argument by combining Ascoli's theorem and the Cantor diagonal process, we get a nondecreasing $\varphi:\N\rightarrow\N$ and
$u\in L^\infty(\R^+,L^2(\R^3))\cap C(\R^+,H^{-3}(\R^3))$ such that: for all $q\in\N$, we have
\begin{equation}\label{theq6}\lim_{n\rightarrow\infty}\|\theta_q(u_{\varphi(n)}-u)\|_{L^\infty([0,T_q],H^{-4})}=0.\end{equation}
\item[$\bullet$] Combining inequalities (\ref{theq3})-(\ref{theq6}) and applying Lemma \ref{lem2.9}, we get
\begin{equation}\label{theq7}u_{\varphi(n)}\rightarrow u\;{\rm strongly\; in}\;L^p_{loc}(\R^+,L^p_{loc}(\R^3)),\;\forall p\in[2,\infty).\end{equation}
\item[$\bullet$] The sequence $(u_n)$ is bounded in the Hilbert space $L^2(\R^+,\dot H^1(\R^3))$, then
    $$u_{\varphi(n)}\rightarrow u\;{\rm weakly\;in}\; L^2(\R^+,\dot H^1(\R^3)).$$
    Particularly $u\in L^2(\R^+,\dot H^1(\R^3))$ and \begin{equation}\label{theq8}\int_0^t\|\nabla u\|_{L^2}^2\leq \liminf_{n\rightarrow\infty}\int_0^t\|\nabla u_{\varphi(n)}\|_{L^2}^2,\;\;\forall t\geq0.\end{equation}
\item[$\bullet$] Prove that for all $T>0$ and $k\in\N$,
\begin{equation}\label{theq10}\int_{0}^T\int_{\R^3}|u|^{2k+2}\leq\liminf_{n\rightarrow\infty}\int_{0}^T\int_{\R^3} |u_{\varphi(n)}|^{2k+2}.\end{equation}
For this take $m\in\N$ and put $\Omega_m=]0,T[\times B_m$. Applying Lemma \ref{lem2.9} to the sequence $({\bf 1}_{\Omega_m}u_{\varphi(n)})$ with $p_0=2k+4$ and $p=2k+2$, with taking into account (\ref{theq4}), we get
$$\int_{0}^T\int_{B_m}|u|^{2k+2}\leq\liminf_{n\rightarrow\infty}\int_{0}^T\int_{B_m} |u_{\varphi(n)}|^{2k+2}.$$
Then
$$\int_{0}^T\int_{B_m}|u|^{2k+2}\leq\liminf_{n\rightarrow\infty}\int_{0}^T\int_{\R^3} |u_{\varphi(n)}|^{2k+2}.$$
By using Monotonic Convergence Theorem, we obtain the desired result.
\item[$\bullet$] Prove that for all $T>0$,
\begin{equation}\label{theq11}\int_0^T\int_{\R^3}(e^{\beta|u|^2}-1)|u|^2\leq \liminf_{n\rightarrow\infty}\int_0^T\int_{\R^3}(e^{\beta|u_{\varphi(n)}|^2}-1)|u_{\varphi(n)}|^2.\end{equation}
For this, take $m\in\N$ and by (\ref{theq10}), we get
$$\int_0^T\int_{\R^3}P_m(\beta|u|^2)|u|^2\leq\liminf_{n\rightarrow\infty}\int_0^T\int_{\R^3}P_m(\beta|u_{\varphi(n)}|^2)|u_{\varphi(n)}|^2.$$
Then$$\int_0^T\int_{\R^3}P_m(\beta|u|^2)|u|^2\leq \liminf_{n\rightarrow\infty}\int_0^T\int_{\R^3}(e^{\beta|u_{\varphi(n)}|^2}-1)|u_{\varphi(n)}|^2.$$
Monotonic Convergence Theorem gives the desired result.
\item[$\bullet$] Combining the above inequalities, we obtain: for all $t\geq0$
\begin{equation}\label{theq12}\|u(t)\|_{L^2}^2+2\int_0^t\|\nabla u(z)\|_{L^2}^2dz+2\alpha\int_0^t\|(e^{\beta |u(z)|^2}-1)|u(z)|^2\|_{L^1}dz\leq \|u^0\|_{L^2}^2.\end{equation}
{\bf It remains to show that $u$ is a solution of the system $(S)$.} To do this, we must verify that
\begin{equation}\label{theq13}\lim_{n\rightarrow\infty}J_{\varphi(n)}\big(u_{\varphi(n)}^iu_{\varphi(n)}^j\big)=u^iu^j\;\;\;{\rm in}\;\;\;{\mathcal D'}(\mathbb [0,\infty)\times\mathbb R^3)),\end{equation}
\begin{equation}\label{theq14}\lim_{n\rightarrow\infty}J_{\varphi(n)}\big[(e^{\beta|u_{\varphi(n)}|^2}-1)u_{\varphi(n)}\big]=
(e^{\beta|u|^2}-1)u\;\;\;{\rm in}\;\;\;{\mathcal D'}(\mathbb [0,\infty)\times\mathbb R^3)).\end{equation}
\item[$\bullet$] Let $\Phi$ in $C^\infty_0(\R^+\times\R^3)$. There is an integer $q\in\N$ such that
$$supp(\Phi)\subset[0,T_q)\times B(0,q).$$
We start by proving (\ref{theq13}). For this, we write
$$\begin{array}{lll}
\Phi u_{\varphi(n)}^iu_{\varphi(n)}^j-\Phi u^iu^j&=&\theta_q\Phi u_{\varphi(n)}^iu_{\varphi(n)}^j-\theta_q\Phi u^iu^j\\
&=&\theta_q(u_{\varphi(n)}^i-u^i)\Phi u_{\varphi(n)}^j+\theta_q(u_{\varphi(n)}^j-u^j)\Phi u^i.
\end{array}$$
Using Lemma \ref{l122}, we get, for $\varepsilon_0\in(0,1)$,
$$\begin{array}{l}
\|\theta_q(u_{\varphi(n)}^i-u^i)\Phi u_{\varphi(n)}^j\|_{L^2(\mathbb R^+,H^{-\varepsilon_0-\frac{1}{2}})}~\quad\quad\quad\quad\quad\quad~\\
~\quad\quad\quad\quad\quad\quad~\leq C\|\theta_q (u_{\varphi(n)}^i-u^i)\|_{L^\infty([0,T_q],H^{-\frac{\varepsilon_0}{2}})}\|\Phi u_{\varphi(n)}^j\|_{L^2([0,T_q],H^{1-\frac{\varepsilon_0}{2}})}\\
~\quad\quad\quad\quad\quad\quad~\leq C\|\theta_q (u_{\varphi(n)}^i-u^i)\|_{L^\infty([0,T_q],H^{-\frac{\varepsilon_0}{2}})}\|\Phi u_{\varphi(n)}^j\|_{L^2([0,T_q],H^{1})}\\
~\quad\quad\quad\quad\quad\quad~\leq C'\|\theta_q(u_{\varphi(n)}^i-u^i)\|_{L^\infty([0,T_q];H^{-\frac{\varepsilon_0}{2}})}\|u_{\varphi(n)}^j\|_{L^2([0,T_q];H^1)}\\
~\quad\quad\quad\quad\quad\quad~\leq C'M_{T_q}\|\theta_q(u_{\varphi(n)}^i-u^i)\|_{L^\infty([0,T_q];H^{-\frac{\varepsilon_0}{2}})},
\end{array}$$
where $M_{T_q}$ is given by (\ref{theq2}). Then, equation (\ref{theq6}) implies $$\displaystyle\lim_{n\rightarrow\infty}\|\theta_q(u_{\varphi(n)}^i-u^i)\Phi u_{\varphi(n)}^j\|_{L^2(\mathbb R^+,H^{-\varepsilon_0-\frac{1}{2}})}=0.$$ In a similar way, we show that
$$\lim_{n\rightarrow\infty}\theta_q(u_{\varphi(n)}^j-u^j)\Phi u^i=0\;{\rm in}\;L^2(\mathbb R^+,H^{-\varepsilon_0-\frac{1}{2}}(\R^3)).$$
Then, we obtain
$$\lim_{n\rightarrow\infty}\Phi u^i_{\varphi(n)}u^j_{\varphi(n)}=\Phi u^iu^j\;\;\;{\rm in}\;\;\;L^2(\mathbb R^+,H^{-\varepsilon_0-\frac{1}{2}}(\R^3)).$$  Using the fact that $(u_n)$ is bounded in $L^\infty(\mathbb R^+,L^2(\R^3))\cap L^2(\mathbb R^+,\dot H^1(\R^3))$, we get
$$\begin{array}{lll}
\|(Id-J_{\varphi(n)})u^i(t)u^j(t)\|_{H^{-\varepsilon_0-\frac{1}{2}}}&\leq&\|(Id-J_{\varphi(n)})u^i(t)u^j(t)\|_{\dot H^{-\varepsilon_0-\frac{1}{2}}}\\
&\leq& (\varphi(n))^{-\varepsilon_0}\|u^i(t)u^j(t)\|_{\dot H^{-\frac{1}{2}}}\\
&\leq& C(\varphi(n))^{-\varepsilon_0}\|u(t)\|_{L^2}\|u(t)\|_{\dot H^1}.
\end{array}$$
Then, we get
$$\|(Id-J_{\varphi(n)})u^i(t)u^j(t)\|_{L^2(\R^+,H^{-\varepsilon_0-\frac{1}{2}})}\leq C(\varphi(n))^{-\varepsilon_0}\|u\|_{L^\infty(\mathbb R^+,L^2)}\|u\|_{L^2(\mathbb R^+,\dot H^1)}.$$
which implies (\ref{theq13}).
\item[$\bullet$] Now, we want to prove (\ref{theq14}). Let $\varepsilon>0$ and $R_0>0$ such that
\begin{equation}\label{eqr}\frac{1}{R_0}<\frac{\varepsilon}{4\|u^0\|_{L^2}^2\|\Phi\|_{L^\infty(\R^+\times\R^3)}+4}.\end{equation}
Let $m_0\in\N$ such that
\begin{equation}\label{eqm}(e^{\beta R_0^2}-1-P_{m_0}(\beta R_0^2))R_0=\sum_{k=m_0+1}^\infty\frac{\beta^k}{k!}R_0^{2k+1}<\frac{\varepsilon}{4\alpha\|\Phi\|_{L^1(\R^+\times\R^3)}+4}.\end{equation}
We have
$$\int_{0}^\infty\int_{\R^3} J_{\varphi(n)}[(e^{\beta |u_{\varphi(n)}|^2}-1)u_{\varphi(n)}]\Phi-\int_{0}^\infty\int_{\R^3} (e^{\beta |u|^2}-1)u\Phi=(\sum_{i=1}^3L_n^i)+S,$$
with
$$\begin{array}{lcl}
L_n^1&=&\displaystyle\int_0^\infty\int_{\R^3} [P_{m_0}(\beta|u_{\varphi(n)}|^2)u_{\varphi(n)}-P_{m_0}(\beta|u|^2)u]\Phi\\
L^2_n&=&\displaystyle\int_0^\infty\int_{\R^3} (Id-J_{\varphi(n)})[(e^{\beta |u_{\varphi(n)}|^2}-1)u_{\varphi(n)}]\Phi\\
L^3_n&=&\displaystyle\int_0^\infty\int_{\R^3} [(e^{\beta |u_{\varphi(n)}|^2}-1-P_{m_0}(\beta|u_{\varphi(n)}|^2))u_{\varphi(n)}]\Phi\\
S&=&\displaystyle\int_0^\infty\int_{\R^3} [(e^{\beta |u|^2}-1-P_{m_0}(\beta|u|^2))u]\Phi.
\end{array}$$
$-$ Estimate of term $L_n^1$: Using Lemma \ref{lempn2}, we get
$$\begin{array}{lcl}|L_n^1|&\leq& \displaystyle C_{m_0,\beta}\int_0^{T_{q}}\int_{B(0,q)} [P_{m_0}(\beta|u_{\varphi(n)}|^2)+P_{m_0}(\beta|u|^2)]\theta_q|u_{\varphi(n)}-u||\Phi|\\
&\leq& \displaystyle C_{m_0,\beta}\|\Phi\|_{L^\infty(\R^+\times\R^3)}\int_0^{T_q}\int_{B(0,q)} [P_{m_0}(\beta|u_{\varphi(n)}|^2)+P_{m_0}(\beta|u|^2)]\theta_q|u_{\varphi(n)}-u|.
\end{array}$$
Applying Cauchy-Schwartz and using Lemma \ref{lempn1}, we get
$$\begin{array}{lcl}
|L_n^1|&\leq& \displaystyle C_{m_0,\beta}\|\Phi\|_{L^\infty(\R^+\times\R^3)}\Big(\|[P_{m_0}(\beta|u_{\varphi(n)}|^2)\|_{L^2([0,T_q]\times B(0,q))}\\
&&+\|P_{m_0}(\beta|u|^2)]\|_{L^2([0,T_q]\times B(0,q))}\Big)\|\theta_q(u_{\varphi(n)}-u)\|_{L^2([0,T_q]\times B(0,q))}\\
&\leq& \displaystyle c_{m_0,\beta}C_{m_0,\beta}\|\Phi\|_{L^\infty(\R^+\times\R^3)}\Big(\|(e^{\beta|u_{\varphi(n)}|^2}-1)|u_{\varphi(n)}|^2\|_{L^1([0,T_q],L^1(\R^3))}^{1/2}\\
&&+\|(e^{\beta|u|^2}-1)|u|^2\|_{L^1([0,T_q],L^1(\R^3))}^{1/2}\Big)\|\theta_q(u_{\varphi(n)}-u)\|_{L^2([0,T_q]\times B(0,q))}.
\end{array}$$
The energy inequalities (\ref{theq1})-(\ref{theq12}) give
$$\begin{array}{lcl}
|L_n^1|&\leq&\displaystyle c_{m_0,\beta}C_{m_0,\beta}\|\Phi\|_{L^\infty(\R^+\times\R^3)}\frac{\sqrt{2}\|u^0\|_{L^2}}{\sqrt{\alpha}}\|\theta_q(u_{\varphi(n)}-u)\|_{L^2([0,T_q],L^2)}.
\end{array}$$
Then, the convergence result (\ref{theq6}) gives an integer $n_1\in\N$ such that: for all $n\geq n_1$, we have
\begin{equation}\label{ineqv1}|L_n^1|<\varepsilon/4.\end{equation}
$-$ Estimate of term $L^2_n$: We have
$$\begin{array}{lcl}|L^2_n|&=& |\displaystyle C\int_0^{T_q}\int_{\R^3}(Id-J_{\varphi(n)})[(e^{\beta |u_{\varphi(n)}|^2}-1)u_{\varphi(n)}]\Phi|\\
&=& |\displaystyle C\int_0^{T_q}((Id-J_{\varphi(n)})[(e^{\beta |u_{\varphi(n)}|^2}-1)u_{\varphi(n)}]|\Phi)_{L^2}|\\
&\leq& \displaystyle C\int_0^{T_q}\|(Id-J_{\varphi(n)})[(e^{\beta |u_{\varphi(n)}|^2}-1)u_{\varphi(n)}]\|_{H^{-4}}\|\Phi(t)\|_{H^4}\\
&\leq& \displaystyle \frac{C}{\varphi(n)}\int_0^{T_q}\|[(e^{\beta |u_{\varphi(n)}|^2}-1)u_{\varphi(n)}]\|_{H^{-3}}\|\Phi\|_{L^\infty(\R^+,H^4)}\\
&\leq& \displaystyle \frac{C}{\varphi(n)}\|[(e^{\beta |u_{\varphi(n)}|^2}-1)u_{\varphi(n)}]\|_{L^1([0,T_q],H^{-3})}\|\Phi\|_{L^\infty(\R^+,H^4)}.
\end{array}$$
By using Proposition \ref{prop1}-(\ref{ps1eq11}) with the choice $(R=1,\;s=3,\;T=T_q)$, we get
$$|L^2_n|\leq\displaystyle \frac{C}{\varphi(n)}\Big(M_{\beta,1}T_q+\frac{1}{2\alpha}\Big)\sigma_{3,3}\|u^0\|_{L^2}^2\|\Phi\|_{L^\infty(\R^+,H^4)},$$
which implies that there is an integer $n_2\in\N$ such that \begin{equation}\label{ineqv2}|L_n^2|<\varepsilon/4,\;\;\forall n\geq n_2.\end{equation}
$-$ Estimate of term $L^3_n$: Put the following set
$$Y_{n,R_0}=\{(t,x)\in[0,T_q]\times B(0,q);\;|u_{\varphi(n)}(t,x)|\leq R_0\}.$$
We have
$$\begin{array}{lcl}
|L^3_n|&\leq&\displaystyle\int_{Y_{n,R_0}} [(e^{\beta |u_{\varphi(n)}|^2}-1-P_{m_0}(\beta|u_{\varphi(n)}|^2))|u_{\varphi(n)}|]|\Phi|\\
&&\displaystyle+\int_{Y_{n,R_0}^c} [(e^{\beta |u_{\varphi(n)}|^2}-1-P_{m_0}(\beta|u_{\varphi(n)}|^2))|u_{\varphi(n)}|]|\Phi|\\
&\leq&\displaystyle [(e^{\beta R_0^2}-1-P_{m_0}(\beta R_0^2))R_0]\|\Phi\|_{L^1(\R^+\times\R^3)}\\
&&\displaystyle+\frac{1}{R_0}\int_{Y_{n,R_0}^c} [(e^{\beta |u_{\varphi(n)}|^2}-1-P_{m_0}(\beta|u_{\varphi(n)}|^2))|u_{\varphi(n)}|^2]\|\Phi\|_{L^\infty(\R^+\times\R^3)}\\
&\leq&\displaystyle [(e^{\beta R_0^2}-1-P_{m_0}(\beta R_0^2))R_0]\|\Phi\|_{L^1(\R^+\times\R^3)}\\
&&\displaystyle+\frac{1}{R_0}\int_{\R^+\times\R^3} [(e^{\beta |u_{\varphi(n)}|^2}-1)|u_{\varphi(n)}|^2]\|\Phi\|_{L^\infty(\R^+\times\R^3)}\\
&\leq&\displaystyle [(e^{\beta R_0^2}-1-P_{m_0}(\beta R_0^2))R_0]\|\Phi\|_{L^1(\R^+\times\R^3)}+\frac{1}{R_0}\|u^0\|_{L^2}^2\|\Phi\|_{L^\infty(\R^+\times\R^3)},
\end{array}$$
Using the choices (\ref{eqr})-(\ref{eqm}), we get
\begin{equation}\label{ineqv3}|L_n^3|<\varepsilon/4,\;\;\forall n\in\N.\end{equation}
$-$ Estimate of term $S$: Put the following subset of $\R^4$,
$$Z_{R_0}=\{(t,x)\in[0,T_q]\times B(0,q);\;|u(t,x)|\leq R_0\}.$$
We have
$$\begin{array}{lcl}
|S|&\leq&\displaystyle\int_{Z_{R_0}} [(e^{\beta |u|^2}-1-P_{m_0}(\beta|u|^2))|u|]|\Phi|\\
&&\displaystyle+\int_{Z_{R_0}^c} [(e^{\beta |u|^2}-1-P_{m_0}(\beta|u|^2))|u|]|\Phi|\\
&\leq&\displaystyle (e^{\beta R_0^2}-1-P_{m_0}(\beta R_0^2))R_0\|\Phi\|_{L^1(\R^+\times\R^3)}\\
&&\displaystyle+\frac{1}{R_0}\Big(\int_{Z_{R_0}^c} (e^{\beta |u|^2}-1-P_{m_0}(\beta|u|^2))|u|^2\Big)\|\Phi\|_{L^\infty(\R^+\times\R^3)}\\
&\leq&\displaystyle (e^{\beta R_0^2}-1-P_{m_0}(\beta R_0^2))R_0\|\Phi\|_{L^1(\R^+\times\R^3)}\\
&&\displaystyle+\frac{1}{R_0}\Big(\int_{\R^+\times\R^3} (e^{\beta |u|^2}-1)|u|^2\Big)\|\Phi\|_{L^\infty(\R^+\times\R^3)}\\
&\leq&\displaystyle (e^{\beta R_0^2}-1-P_{m_0}(\beta R_0^2))R_0\|\Phi\|_{L^1(\R^+\times\R^3)}+\frac{1}{2\alpha R_0}\|u^0\|_{L^2}^2\|\Phi\|_{L^\infty(\R^+\times\R^3)}.
\end{array}$$
By the choices (\ref{eqr})-(\ref{eqm}) we get
\begin{equation}\label{ineqv4} |S|<\varepsilon/4.\end{equation}
Combining (\ref{ineqv1}),\,(\ref{ineqv2}),\,(\ref{ineqv3}) and (\ref{ineqv4}), we get
$$\Big|\int_{0}^\infty\int_{\R^3} J_{\varphi(n)}[(e^{\beta |u_{\varphi(n)}|^2}-1)u_{\varphi(n)}]\Phi-\int_{0}^\infty\int_{\R^3} (e^{\beta |u|^2}-1)u\Phi\Big|<\varepsilon,\;\forall n\geq\max(n_1,n_2),$$
which implies (\ref{theq14}) and the proof of Theorem \ref{th1} is finished.
\end{enumerate}

\end{document}